\documentclass[12pt]{amsart}

\makeatletter \@addtoreset{equation}{section}
\@addtoreset{figure}{section}

\makeatother

\usepackage{xypic,amscd}
\usepackage{amsmath}
\usepackage{amssymb}
\usepackage{latexsym}
\newcommand{\Pic}{\operatorname{Pic}}

\newcommand{\mult}{\operatorname{mult}}
\newcommand{\Sing}{\operatorname{Sing}}

\newcommand{\Cr}{\operatorname{Cr}}
\newcommand{\Bs}{\operatorname{Bs}}
\newcommand{\PSL}{\operatorname{L}}
\newcommand{\Aut}{\operatorname{Aut}}
\newcommand{\di}{\operatorname{d}}

\newcommand{\OOO}{{\mathcal O}}
\newcommand{\PP}{\mathbb{P}}

\newcommand{\FF}{\mathbb{F}}

\newcommand{\CC}{\mathbb{C}}
\newcommand{\QQ}{\mathbb{Q}}
\newcommand{\AAA}{\mathfrak{A}}
\newcommand{\SSS}{\mathfrak{S}}

\newcommand{\comment}[1]{}

\newtheorem{theorem}[equation]{Theorem}
\newtheorem{proposition}[equation]{Proposition}
\newtheorem{lemma}[equation]{Lemma}
\newtheorem{corollary}[equation]{Corollary}
\newtheorem{claim}[equation]{Claim}

\theoremstyle{definition}
\newtheorem{example}[equation]{Example}
\newtheorem{definition}[equation]{Definition}
\newtheorem{notation}[equation]{Notation}
\newtheorem{construction}[equation]{Construction}

\newtheorem{remark}[equation]{Remark}

\newtheorem{case}{}[equation]

\theoremstyle{remark}

\date{}

\begin{document}

\title[]{Log del Pezzo surfaces with simple automorphism groups}
\author[]{Grigory Belousov}
\thanks{The work was partially supported by grant N.Sh.-1987.2008.1}
\maketitle

\section{introduction}

A (singular) \emph{del Pezzo surface} is a normal projective surface
$X$ over $\CC$ whose anticanonical divisor $-K_{X}$ is ample.

In this paper we consider del Pezzo surfaces with only log terminal
singularities admitting an action of a finite simple group $G$. Any
del Pezzo surface with log terminal singularities is rational (see,
e.g., \cite{AN}). Hence such a group $G$ is contained in $\Cr(2)$,
the plane Cremona group. Finite subgroups of $\Cr(2)$ are classified
(see \cite{DI}). By \cite{DI}, There are only three finite simple
subgroups of $\Cr(2)$: $\AAA_{5}$, $\AAA_{6}$ and $\PSL_{2}(7)$,
where $\PSL_{2}(7)$ is the simple group of order $168$ and
$\AAA_{n}$ is the alternating group. These groups have their names:
$\PSL_{2}(7)$ is called Klein's simple group, $\AAA_{6}$ is called
the Valentiner group, and $\AAA_{5}$ is the icosahedral one. In the
present paper we classify del Pezzo surfaces admitting an action of
one of these groups.

\begin{example}[see, e.g. {\cite{DI}}]
\label{Ex0} The Klein group $G=\PSL_{2}(7)$ has an irreducible
three-dimensional representation, so $G$ acts on the projective
plane. There is an invariant quartic, so-called \textit{Klein
quartic} $C^{\mathrm{k}}:=\{x^{3}y+y^{3}z+z^{3}x=0\}$. Consider the
double cover $S_2^{\mathrm{k}}\to \PP^2$ ramified along
$C^{\mathrm{k}}$. Then $S_2^{\mathrm{k}}$ is a smooth del Pezzo
surface of degree two. The action of $G$ lifts naturally to
$S_2^{\mathrm{k}}$. Such $S_2^{\mathrm{k}}$ is unique up to
isomorphism and $\rho(S_2^{\mathrm{k}})^G=1$.
\end{example}

\begin{example}
\label{Ex1} Let $n$ be a positive integer and let
$k\in\{12,20,30,60\}$. The group $G=\AAA_{5}$ acts naturally on the
Hirzebruch surface $\FF_{2n}$. Let $M$ and $D$ be disjoint sections
with $M^2=2n$, $D^2=-2n$. Let $p_{1}:Z_{1}\rightarrow\FF_{2n}$ be
the blowup of an orbit consisting of $k$ points on $M$. Let
$p_{2}:Z_{2}\rightarrow Z_{1}$ be the blowup of an orbit consisting
of $k$ points on proper transform of $M$ and so on. Let
$r:Z_{a}\rightarrow F_{2n,ak-2n,a}$ be the contraction of the curves
with self-intersection number is less than $-1$, where $a$ is the
number of blowups of points on proper transform of $M$ and
$ak-2n>0$. Then $F_{2n,ak-2n,a}$ is a del Pezzo surface admitting a
non-trivial action of $G$ and $\rho(F_{2n,ak-2n,a})^G=1$. Note that
$F_{2n,ak-2n,a}$ has two singular points of types
$\frac{1}{2n}(1,1)$, $\frac{1}{ak-2n}(1,1)$ and $k$ Du Val
singularities of type $A_{a-1}$. It is possible that $a=1$ and then
the singular locus of $F_{2n,k-2n,1}$ consists of two points.
\end{example}

\begin{example}
\label{Ex2} Let $k\in\{12,20,30,60\}$. The group $G=\AAA_{5}$ acts
naturally on $\PP^2$. Let $C$ be (a unique) $G$-invariant conic on
$\PP^2$. Let $p_{0}:Z_{0}\rightarrow\PP^2$ be the blowup of an orbit
consisting of $k$ points $P_{1},\dots, P_{k}$ on $C$ and let $C_{0}$
be the proper transform of $C$. Let $p_{1}:Z_{1}\rightarrow Z_{0}$
be the blowup an orbit of points on $C_{0}$ that correspond to
$P_{1},\dots, P_{k}$. Repeating this procedure $s+1$ times we obtain
a smooth surface $Z_{s}$. Let
$r:Z_{s}\rightarrow\tilde{\PP^2}_{k,s}$ be the contraction of all
rational curves whose self-intersection number is at most $-2$. Then
$\tilde{\PP^2}_{k,s}$ is a del Pezzo surface admitting a non-trivial
action of $G$ and $\rho(\tilde{\PP^2}_{k,s})^G=1$. The singular
locus of $\tilde{\PP^2}_{k,s}$ consists of one (fixed) point of type
$\frac{1}{k(s+1)-4}(1,1)$ and $k$ Du Val singular points of type
$A_{s}$. It is possible that  $s=0$ and then the singular
locus of $\tilde{\PP^2}_{k,0}$ consists of one (fixed) point.
\end{example}

The main result of this paper is the following:

\begin{theorem}
\label{1.1} Let $X$ be a del Pezzo surface with log terminal
singularities and let $G\subset \Aut(X)$ be a finite simple group.
\begin{enumerate}
\item If $G\simeq \AAA_{5}$ and $\rho(X)^G=1$, then there are the following cases:
\begin{itemize}
\item $X\simeq\PP^2$.
\item $X\simeq S_{5}$, where $S_{5}$ is a smooth del Pezzo surface of
degree $5$.
\item $X\simeq\PP(1,1,2n)$, a cone over a rational normal curve of degree $2n$.
\item $X\simeq F_{2n,ak-2n,a}$, see Example \ref{Ex1}.
\item $X\simeq\tilde{\PP^2}_{k,s}$, see Example \ref{Ex2}.
\end{itemize}
\item
If $G$ is the Klein group, then $X\simeq\PP^2$ or $X\simeq S_2^{\mathrm{k}}$.
\item
If $G$ is the Valentiner group, then $X\simeq\PP^2$.
\end{enumerate}
\end{theorem}


The author is grateful to Professor Y. G. Prokhorov for suggesting
me this problem and for his help. The author would like to thank
Professor I. A. Cheltsov for useful comments.

\section{Preliminaries}
\begin{notation}
\label{not} We work over $\CC$. Throughout this paper $G$ is one of
the following groups: $\AAA_{5}$, $\AAA_{6}$, or $\PSL_{2}(7)$. $X$
denotes a del Pezzo surface with at worst log terminal singularities
admitting a non-trivial action of $G$. We also employ the following
notation:
\begin{itemize}
\item $\FF_{n}$: the Hirzebruch surface, $\FF_{0}\simeq \PP^1\times\PP^1$.
\item $\PP(a,b,c)$: weighted projective plane.
\item $S_d$: a del Pezzo surface of degree $d$.
\item $\rho(X)$: the Picard number.
\item
$G$-surface: a surface $V$ with a given embedding $G\subset \Aut (V)$.
\item $\rho(X)^G$: the $G$-invariant Picard number.
\item $(n)$-curve: a smooth rational curve whose self-intersection number equals to $n$.
\end{itemize}
\end{notation}

\begin{definition}
\label{DCoD} Let $S$ be a normal projective surface and let
$f:\tilde{S}\rightarrow S$ be a resolution. Let $D=\sum_{i}D_{i}$ be
the exceptional divisor. Then there is a unique $\QQ$-divisor
$D^\sharp=\sum_{i}\alpha_{i}D_{i}$ such that $f^{*}K_{S}\equiv
K_{\tilde{S}}+D^\sharp$. The numbers $\alpha_{i}$ are called
\emph{codiscrepancy} of $D_{i}$.
\end{definition}

\begin{lemma}
\label{fixed-point}
Let $V$ be a $G$-surface with at worst log terminal
singularities and let $P\in V$ be a fixed point. Then $P$ is
singular and $G\simeq\AAA_{5}$. Moreover, $P$ has type
$\frac{1}{r}(1,1)$, where $r$ is even.
\end{lemma}

\begin{proof}
Let $P\in V$ be a $G$-fixed point. Assume that $P$ is a smooth
point. Then $G$ acts on the Zariski tangent space $T_{P,V}$. Since
$G$ is a finite simple group, we see that $G$ has no non-trivial
two-dimensional representations. Hence $P$ is singular. Let
$\tilde{V}\rightarrow V$ be the minimal resolution of $P$ and let
$D=\sum D_{i}$ be the exceptional divisor. Then $G$ acts on $D$.
Since $G$ does not admit any embeddings to $\SSS_{k}$, where $k\leq
4$, we see that $D$ consists of one irreducible component. Hence $P$
has type $\frac{1}{r}(1,1)$. On the other hand, the Klein group and
the Valentiner group do not admit a non-trivial action on a smooth
rational curve. Hence $G\simeq\AAA_{5}$.

Finally, the action of $\AAA_5$ on $\tilde{V}$
induces an action of $\AAA_5$ on the total space of the  conormal bundle
$N_{D/\tilde{V}}^{\vee}\simeq \OOO_{\PP^1}(r)$. In particular, the group
$\AAA_5$ naturally acts on
\[
H^0(D, N_{D/\tilde{V}}^{\vee})\simeq H^0(\PP^1,\OOO_{\PP^1}(r))
\simeq S^r H^0(\PP^1,\OOO_{\PP^1}(1)).
\]
This is possible if and only if $r$ is even.
\end{proof}

\begin{lemma}
\label{CB}
Let $G\simeq\AAA_{5}$ and let $f:S\rightarrow\PP^1$ be a
smooth relatively minimal conic bundle with an action of $G$. Then
$S\simeq\FF_{2n}$. Moreover, there are two possibilities:
\begin{itemize}
\item
$S=\PP^1\times\PP^1$ with non-trivial action on each factor;
\item
$S=\FF_{2n}$, $n\geq 0$, there is an invariant section and this case occurs for every $n\geq 0$.
\end{itemize}
\end{lemma}

\begin{proof}
Let $$\alpha:G\rightarrow O(\Pic(S)),\quad g\rightarrow
(g^{*})^{-1}$$ be the natural representation of $G$ in the
orthogonal group of $\Pic(S)$. By \cite[Theorem 5.7]{DI}, we have
$\ker(\alpha)\neq\{e\}$. Since $G$ is a simple group, we see that
$\ker(\alpha)=G$. Hence $f$ has no singular fibers. Then
$S=\FF_{r}$.

Since the case $r=0$ is
trivial, we assume that $r>0$. Consider the contraction $\varphi:
\FF_r\to V_r$ of the negative section. Here $V_r$ is a cone in
$\PP^{r+1}$ over a rational normal curve $C_r\subset \PP^r$ of
degree $r$ or, equivalently, the weighted projective plane
$\PP(1,1,r)$. Clearly, $\varphi$ is $\AAA_5$-equivariant, so
$\AAA_5$ acts non-trivially on $V_r$.
By Lemma \ref{fixed-point} $r$ is even.
On the
other hand, one can write down an action of $\AAA_5$ on
$\PP(1,1,2n)$ explicitly.
\end{proof}

\begin{lemma}
\label{Discr}
Let $P\in X$ be a log terminal singularity and let
$f:\tilde{X}\rightarrow X$ be its minimal resolution. Let
$\sum\alpha_{i}D_{i}$ be a $\QQ$-divisor such that
$$
f^{*}K_{X}\equiv K_{\tilde{X}}+\sum\alpha_{i}D_{i}.
$$
Assume that
$\alpha_{i}<\frac{2}{5}$ for every $i$. Then $P$ is either a Du Val
singularity or $P$ has type $\frac{1}{3}(1,1)$, i.e. the exceptional
divisor of $f$ consists of a single $(-3)$-curve.
\end{lemma}

\begin{proof}
Assume that there is a component $D_{j}$ of $\sum D_{i}$ such that
$D_{j}^2\leq -4$. Then by {\cite[Lemma 2.17]{Al}} we have
$\alpha_{j}\geq\frac{1}{2}$, a contradiction.

Hence, $D_{i}^2\geq -3$
for every $i$. Assume that there is a component $D_{j}$ with
$D_{j}^2=-3$ and a component $D_{k}$ with $D_{j}\cdot D_{k}=1$. Then
by {\cite[Lemma 2.17]{Al}} we have $\alpha_{j}\geq\frac{2}{5}$.
Again we have a
contradiction.

Therefore, $P$ is either a Du Val or
the exceptional divisor $\sum D_{i}$ has only one component $D_1$ with
$D_1^2=-3$.
\end{proof}

The following Lemma is a consequence of the classification of log
terminal singularities (see \cite{Br}).

\begin{lemma}
\label{Discr2} Let $X$ be a projective normal surface. Let $P\in X$
be a log terminal non-Du Val singularity and let
$f:\tilde{X}\rightarrow X$ be its minimal resolution. Let
$\sum\alpha_{i}D_{i}$ be a codiscrepancy $\QQ$-divisor over $P$.
Assume that there is a $(-1)$-curve $E$ and a morphism
$g:\tilde{X}\rightarrow Z$ such that $g(E+\sum D_{i})$ is a smooth
point. Then $E\cdot\sum\alpha_{i}D_{i}\ge\frac{2}{11}$.
\end{lemma}

\begin{proof}
Consider minimal resolution of log terminal singularities
{\cite{Br}} case by case. For example if $P\in X$ has type
$\frac{1}{7}(1,3)$, then $E\cdot\sum\alpha_{i}D_{i}=\frac{2}{7}$.
\end{proof}

\begin{proposition}
\label{Pr1} In notation \ref{not} assume that $X$ has at worst Du
Val singularities. Then we have one of the following cases:
\begin{enumerate}
\item $G$ is the Valentiner group and $X\simeq\PP^2$.
\item $G$ is the Klein group and $X\simeq\PP^2$ or $X\simeq S_2^{\mathrm{k}}$.
\item
$G\simeq \AAA_{5}$ and either $X$ is smooth or $X\simeq\PP(1,1,2)$.
If moreover $\rho(X)^G=1$, then
$X$ is isomorphic to $\PP^2$, $S_{5}$ or $\PP(1,1,2)$.
\end{enumerate}
\end{proposition}

\begin{proof}
We use some elementary facts on del Pezzo surfaces with Du Val singularities
(see e.g. \cite{HW}).
Recall that $9\geq K_{X}^2\geq 1$ and
$$\dim H^{0}(X,\OOO(-K_{X}))=K_{X}^2+1.$$
So, we have the following
cases:

\begin{case}
$K_{X}^2=1$. In this case $\dim |-K_{X}|=1$ and there is a
non-singular fixed point $\{p\}=\Bs|-K_{X}|$. On the other hand, by
Lemma \ref{fixed-point} every fixed point is singular, a
contradiction.
\end{case}

\begin{case}
$K_{X}^2=2$. In this case $\dim |-K_{X}|=2$. The linear system
$|-K_{X}|$ defines a double cover $\phi=\phi_{|-K_{X}|}:
X\rightarrow\PP^2$. Let $B\subset\PP^2$ be the ramification divisor
of $\phi$. We have $\deg B=4$. Since $G$ is simple, we see that $B$
is irreducible. Since the number of singular points of $B$ is at
most three, we see that $B$ is smooth. So is $X$. By \cite{DI} we
have $G\simeq \PSL_{2}(7)$ and $X\simeq S_2^{\mathrm{k}}$.
\end{case}

\begin{case}
$K_{X}^2=3$. In this case $\dim |-K_{X}|=3$ and $X=X_{3}\subset\PP^3$
is a cubic surface. Here $G$ has a faithful representation in
$H^0(X,-K_{X})=\CC^4$. Assume that $G$ is the Klein group or the
Valentiner group. Then $G$ has no irreducible four-dimensional
representations. So, the representation on $H^0(X,-K_{X})$ is
reducible. Hence there is a $G$-invariant hyperplane $H$.
The intersection $H\cap X$ is a ($G$-invariant) smooth elliptic curve because otherwise
we get a fixed point $P\in H\cap X$ which is impossible.
Since a simple group cannot act on an elliptic curve, we get a
contradiction. Hence $G\simeq\AAA_{5}$.

We claim that the natural representation of $G$ on
$H^0(X,-K_{X})=\CC^4$ is irreducible. Indeed, otherwise there is an
invariant hyperplane $H\subset \PP^3$ and as above we get a fixed
point $P\in H\cap X$. Consider the representation of $G$ on the
Zariski tangents space $T_{P,X}$. Since $\AAA_{5}$ has no
irreducible two-dimensional representations, $\dim T_{P,X}=3$, i.e.
the point $P\in X$ is singular (and Du Val). Take a $G$-equivariant
local embedding $(X,P) \hookrightarrow \CC^3=T_{P,X}$ into the
corresponding affine chart. Let $f=f_{2}+f_{3}$ be the local
equation of $X$ at $P$, where $f_{i}$ is a homogeneous polynomial of
degree $i$. We have $f_{i}\neq 0$ and $\AAA_{5}$ acts on $\CC^3$ so
that $f_{i}$ are invariants. Therefore, $\AAA_{5}$ acts on $\PP^2$
so that the locus $\{f_{2}=0\}$ is an invariant conic and
$\{f_{2}=0\}\cap\{f_{3}=0\}$ is an invariant subset consisting of
$\le 6$ points. On the other hand, any orbit of $\AAA_5$ on a smooth
rational curve contains at least $6$ points. The contradiction
proves our claim. Hence there are no fixed points of $G$ on $X$.

Thus the representation of $G$ on $H^0(X,-K_{X})=\CC^4$
is irreducible.
This representation can be regarded as an invariant hyperplane
$\sum x_i=0$ in $\CC^5$, where $\AAA_5$ acts on $\CC^5$ by
permutations of coordinates.
The ring of invariants $\CC[x_1,\dots,x_5]^{\AAA_{5}}$ is generated by
$\sigma_{1}$, \dots, $\sigma_{5}$, $\delta$,
where $\sigma_{i}$ is the symmetric polynomial of degree $i$ and $\delta$ is
the discriminant.
Therefore, the equation of our cubic surface $X\subset \PP^3\subset \PP^4$ can be written as
$\sigma_{3}=\sigma_{1}^3=0$. This surface is smooth.
\end{case}

\begin{case}
$K_{X}^2\geq 4$. In this case the linear system $|-K_X|$ defines an embedding
$X\hookrightarrow\PP^{K_{X}^2}$. Let $\pi:X_{0}\rightarrow X$ be the minimal
resolution.
Then $\rho(X_0)=10-K_{X_0}^2=10-K_{X}^2\le 6$.
So the number of singular points of $X$ is at most $5$.
Moreover, if $X$ has exactly $5$ singular points, then
$K_{X}^2=4$ and $\rho(X)=1$.
On the other hand, such a surface $X$ does not exist
(see e.g. \cite{F}, \cite{B}). Hence $X$ has at most $4$ singular points.

Assume that $G$ is the Klein group or the Valentiner
group. Run the $G$-equivariant MMP on $X_{0}$. At the end we obtain a del Pezzo
surface $S$ with $\rho(S)^G=1$ and $K_S^2\ge 4$.
So, $S\simeq \PP^2$ \cite{DI}.
Let $s:=\rho(X/X_0)$, the number of exceptional curves of $X\to S$.
Since the Klein group and the Valentiner group do not admit any embeddings to
$\SSS_{5}$ and do not act non-trivially on a rational curve, we see that $s\geq 6$. Therefore,
\[
 K_{X}^2=K_{X_{0}}^2=K_{S}^2-s\leq 9-s<4,
\]
a contradiction.

Thus, $G\simeq\AAA_{5}$. Assume that $X$ is singular and $X\not
\simeq\PP(1,1,2)$. Since $X$ has at most $4$ singular points, there
is a singular fixed point $P$ of $G$ on $X$. Note that there is a
line $E_{1}\subset X\subset\PP^d$ passing through $P$, an image of a
$(-1)$-curve $\ell\subset X_0$. Therefore, there is an orbit of
lines $E_{1},\dots, E_{p}$  passing through $P$, where $p\geq
5$. On the other hand, $X\subset \PP^{K_{X}^2}$ is an intersection
of quadrics (see e.g. \cite{HW}). Hence there are at most four lines
on $X\subset\PP^d$ passing through $P$, a contradiction.
\end{case}
The last assertion follows by \cite{DI}.
\end{proof}

\begin{definition}
Let $S$ be a normal projective surface and let $\Delta$ be an
effective $\QQ$-divisor on $S$. We say that $(S,\Delta)$ is a
\textit{weak log del Pezzo surface} if the pair $(S,\Delta)$ is
Kawamata log terminal (klt) and the divisor $-(K_S+\Delta)$ is nef
and big.
\end{definition}

\begin{remark}
\begin{enumerate}
 \item
Let $(S,\Delta)$ be a weak log del Pezzo surface and let
$\varphi: S\to S'$ be a birational contraction to a normal surface $S'$.
Then $(S',\varphi_*\Delta)$ is also a weak log del Pezzo surface.
 \item
For any weak log del Pezzo surface $(S,\Delta)$ the Mori cone
$\operatorname{NE}(S)$ is polyhedral and generated by contractible
extremal rays.
\end{enumerate}
\end{remark}

\begin{construction}
\label{constr1} Under notation \ref{not}, let $\pi:X_{0}\rightarrow
X$ be the minimal resolution. Run the $G$-equivariant MMP on
$X_{0}$. We obtain a sequence of birational contractions of smooth
surfaces $\phi_{i}: X_i\to X_{i+1}$. At the last step $\phi_{p}:
X_p\to X_{p+1}$ we have either a conic bundle over $X_{p+1}\simeq
\PP^1$ or a contraction to a del Pezzo surface $X_{p+1}$ with
$\rho(X_{p+1})^G=1$ (see \cite{DI}). Write $\pi^{*}K_{X}\equiv
K_{X_{0}}+\Delta_{0}$, where $\Delta_{0}$ is an effective
$\pi$-exceptional $\QQ$-divisor. Note that $(X_0,\Delta_0)$ is a
weak log del Pezzo surface. Define by induction
$\Delta_{i}=\phi_{i-1*}\Delta_{i-1}$. On each step of the MMP the
above property is preserved: $(X_i,\Delta_i)$ is is also a weak log
del Pezzo surface. Since $\rho(X_{p})^G=2$, we see that there is a
$G$-equivariant extremal contraction $g: X_{p}\rightarrow Y$ such
that $g$ is different from $\phi_{p}$. Thus we get the following
sequence of $G$-equivariant contractions:
$$
\begin{CD}
X_{0}@>\phi_{0}>>X_{1}@>\phi_{1}>>\dots @>\phi_{p-1}>> X_{p}@>\phi_{p}>> X_{p+1}
\\
@. @.@. @V{g}VV
\\
@. @. @. Y
\end{CD}
$$
We distinguish the following cases:
\begin{enumerate}
\item
$Y$ is a curve. Then $Y\simeq\PP^1$ and $g$ is a conic bundle
with $\rho(X_{p}/Y)^G=1$. Moreover, in this case $X_p$ is a smooth
del Pezzo surface with $\rho(X_p)^G=2$.
Since the groups $\AAA_6$ and $\PSL_2(7)$ cannot act non-trivially on a rational curve,
we have $G\simeq \AAA_5$.
\item
$Y$ is a smooth surface. Then the contraction $g$ is $K$-negative.
In this case both $Y$ and $X_p$ are smooth del Pezzo surfaces
with $\rho(X_p)^G=2$ and $\rho(Y)^G=1$.
By Proposition \ref{Pr1} we have $G\simeq \AAA_5$.
\item
$Y$ is a singular surface. Then $(Y,g_*\Delta_p)$ is a weak log del
Pezzo surface. In particular, $Y$ is a del Pezzo surface with log
terminal singularities. The group $G$ transitively acts on
$\Sing(Y)$.
\end{enumerate}
Assume that both contractions $g$ and $\phi_p$ are birational. Let
$D=\sum_{i=1}^m D_i$ be the $g$-exceptional divisor, let
$B_i=\phi_p(D_i)$, and let $B=\sum_{i=1}^m B_i$. Since
$\rho(X_{p})^G=2$, we see that the group $G$ acts transitively on
$\{D_i\}$ and on $\{B_i\}$, so the curves $B_i$ have the same
anti-canonical degrees and self-intersection numbers. Since
$\rho(X_{p+1})^G=1$, the divisor $B$ is ample and proportional to
$-K_{X_{p+1}}$. Hence $B$ is connected. Assume that
$\Sing(B)=\emptyset$. Then $B$ is an irreducible curve. Since $B$ is
rational, by the genus formula $K_{X_{p+1}}+B$ is negative. This is
possible only if $X_{p+1}\simeq \PP^2$. Thus we have the following.
\end{construction}

\begin{claim}
In the above notation either
\begin{enumerate}
\item
$\Sing(B)\neq\emptyset$, or
\item
$X_{p+1}\simeq \PP^2$ and $B$ is a smooth irreducible curve of degree $\le 2$.
\end{enumerate}
\end{claim}

\begin{construction}
\label{constr3} Under the notation of \ref{not}, let $\rho(X)^G=1$.
Assume that $X$ is singular. Consider the minimal resolution $\mu
:Y\to X$ and let $R=\sum_{i=1}^n R_i$ be the exceptional divisor.
The action of $G$ lifts naturally to $Y$. Write
\[
\mu^*K_X=K_Y+\sum\alpha_i R_i,
\]
where $0\le \alpha_i<1$. Fix a component, say
$R_1\subset R$, and let $R'=R_1+\cdots+R_k$ be its $G$-orbit. We can
contract all the curves in $R-R'$ over $X$:
\[
\mu : Y\stackrel{\eta}{ \longrightarrow} \bar X \stackrel{\phi}{\longrightarrow} X.
\]
Then $\rho(\bar X)^G=2$ and
\[
\phi^*K_X=\eta_*K_Y=K_{\bar X}+ \Delta,\qquad
\text{where $\Delta:=\eta_*\sum \alpha_iR_i$.}
\]
Therefore, $(\bar X, \Delta)$ is a weak log del Pezzo surface.
Let $\psi : \bar X\to X'$ be (a unique) $K_{\bar X}+ \Delta$-negative contraction.
Clearly, $\psi\neq \phi$, $\psi$ does not contract any component of $\Delta$, and
$\psi$ is also $K_{\bar X}$-negative.
We get the following $G$-equivariant diagram:
\[
\xymatrix{
&\bar X\ar[dl]_{\phi}\ar[dr]^{\psi}&
\\
X&&X', }
\]
where $X'$ is either a smooth rational curve or
a del Pezzo with at worst log terminal singularities and $\rho(X')^G=1$.

For a normal surface $V$, denote by $\di(V)$ the Picard number of
its minimal resolution. In our situation, $\psi \circ \eta$ is a
non-minimal resolution of singularities (because $-K_{\bar X}$ is
$\psi$-ample). Hence $\di(X')<\di (X)$.
\end{construction}

The following procedure is well-known.
It is called the ``2-ray game''.

\begin{construction}
\label{Const} Apply our construction \ref{constr3} several times.
We get the following sequence of $G$-equivariant birational
morphisms:
$$
\xymatrix{& \bar{X}_{0}\ar[dl]_{\phi_0} \ar[dr]^{\psi_0}&
&& &\bar{X}_{d}\ar[dl]_{\phi_{d}} \ar[dr]^{\psi_{d}} \\
**[r]X=X_{0} && X_{1}&\dots & X_{d} && X_{d+1},}
$$
Since $\di (X_i)>\di (X_{i+1})$, the process terminates. Thus at the
end we get $X_{d+1}$ which is either a smooth curve or a smooth del
Pezzo surface with $\rho(X_{d+1})^G=1$. Recall that each $X_i$ for
$i=1,\dots,d$ is a del Pezzo surface with log terminal singularities
and $\rho(X_{i})^G=1$.

Note that on each step the extraction $\phi_{i}$ is not unique; this
obviously depends on the choice $R'$ (in notation of \ref{constr3}).
For our purposes it is convenient to choose $R'$ in one of the following
ways:

\begin{case}
\label{2.11.1} $R'$ is the orbit of exceptional curves over
non-fixed points with maximal codiscrepancy.
\end{case}

\begin{case}
\label{2.11.2} $\phi_{i}(R')$ is a fixed point $P\in X_{i}$. By
Lemma \ref{fixed-point} $R'$ is a unique exceptional curve over $P$.
\end{case}
\end{construction}

\section{The Valentiner and Klein groups}
In this section we prove our main theorem in the case, where $G=
\AAA_{6}$ or $\PSL_2(7)$ (i.e., $G$ is the Valentiner or Klein
group).

\begin{proposition}
\label{1/3} Assume that the surface $X$ is singular and $G$ is
either the Klein group or the Valentiner group. Then $X_p$ has only
cyclic quotient singularities of type $\frac13(1,1)$.
\end{proposition}

\begin{proof}
Apply construction \ref{constr1}. By our assumption we get the case
(iii), i.e., the contraction $g$ is birational and $Y$ is a singular
del Pezzo surface (with log terminal singularities and
$\rho(Y)^G=1$). Moreover, the contraction $\phi_{p+1}$ is also
birational and the exceptional loci of $g$ and $\phi_p$ are
reducible (because $G$ cannot act non-trivially on a rational
curve). Write
\[
g^{*} K_{Y}=K_{X_{p}}+\sum_{i=1}^m\alpha_{i}D_{i}\quad
\text{and}\quad D=\sum_{i=1}^m D_{i},
\]
where, as above, $D_{i}$'s are $g$-exceptional curves. Since the
group $G$ acts transitively on $\{D_i\}$, we have
\[
\alpha_{1}=\dots=\alpha_{m}:=\alpha\quad \text{and}\quad -D_1^2=\dots=-D_m^2:=n.
\]
Further, by the classification of log terminal singularities
\cite{Br} the exceptional divisor over every singular point is
either a pair of $(-n)$-curves or a single $(-n)$-curve (otherwise
$G$ cannot interchange the $D_i$'s).

We claim that $n\leq 3$. Indeed, assume that $n>3$. Note that
$\phi_p$ is the blowup of points in $\Sing(B)$. Let $E$ be a
$\phi_p$-exceptional curve on $X_{p}$. Then
$$
0>K_{Y}\cdot g_{*}E= g^{*}K_{Y}\cdot
E=\left(K_{X_{p}}+\sum_{i=1}^m\alpha D_{i}\right)\cdot
E\geq-1+2\alpha.
$$
Therefore, $\alpha<\frac{1}{2}$. On the other hand,
$$
0=g^{*}K_{Y}\cdot D_j=\left(K_{X_{p}}+\sum_{i=1}^m\alpha
D_{i}\right)\cdot D_j\ge n-2 -\alpha n, \quad n\alpha\ge n-2.
$$
Hence, $n\leq 3$. Moreover, if $n=3$, then $\alpha D_j\cdot \sum
D_i=-1$ and so $D_j\cdot \sum_{i\neq j} D_i<1$. This means that
$D_i\cap D_j=0$ for $i\neq j$. Hence, every singular point on $Y$ is
either Du Val of type $A_l$, $l\le 2$, or a cyclic quotient
singularity of type $\frac{1}{3}(1,1)$. By Proposition \ref{Pr1} we
are done.
\end{proof}

By Proposition \ref{Pr1} we may assume that the singularities of $X$
are worse that Du Val. Apply construction \ref{constr1}. We get the
case (iii). In particular, $\Sing(Y)\neq \emptyset$ and is a del
Pezzo surface with log terminal singularities and $\rho(Y)^G=1$.
Moreover, $X_{p+1}\simeq\PP^2$ or $S_2^{\mathrm k}$ and the latter
is possible only for $G= \PSL_2(7)$ (see \cite{DI}). As in the proof
of Proposition \ref{1/3}, let $D=\sum D_i$ is the $g$-exceptional
divisor and let $B_i:=\phi_p(D_i)$. By Proposition \ref{1/3} every
singular point on $Y$ is of type $\frac{1}{3}(1,1)$, i.e. $D$
consists of disjoint $(-3)$-curves.

First we consider the case $X\simeq\PP^2$.
Then $B_i\cap B_j\neq \emptyset$ and so
$\phi_p$ is a blowup of points in $B_i\cap B_j$, $i\neq j$.
We claim that every curve $B_i$ is smooth and
there are at most two components of $B$ passing through
every point $P\in \PP^2$.
Indeed, assume the
converse. Then
$$
0>g^{*} K_{Y}\cdot E=(K_{X_{p}}+\textstyle \frac{1}{3}\sum D_{i})\cdot
E\geq-1+3\cdot\frac{1}{3}= 0.
$$
Therefore, every $B_i$ is smooth.
Further, since the curves $B_i$ are rational, $k\le 2$.
If $k=1$, then the $B_i$ are lines and on every line we blow up four
points. Hence, number of these lines is equal to five, a contradiction.

Finally, consider the case $k=2$. Then the $B_i$ are smooth conics and
on every conic we blow up seven
points. It is easy to see that the number of points of intersection of conics
is divisible by four, a contradiction.

Now consider the case
$X_{p+1}\simeq S_2^{\mathrm{k}}$.
Then $G$ is the Klein group.
Let $r:=\rho(X_p/X_{p+1})$.
Recall that $m$ is
the number of $D_i$'s.
Then by Noether's formula
$$
0<K_{Y}^2=g^{*} K_{Y}^2=K_{X_p}^2 -D^3
= 10-\rho(X_{p})+\textstyle \frac{m}{3}=2-r+\frac{m}{3}.
$$
Since $m\leq r+7=\rho(X_p)-1$, we see that
\[
0<2-r+\frac{m}{3}\leq-\frac{2m+1}{3}<0,
\]
a contradiction.

\section{The icosahedral group}
It remains to consider the case $G\simeq \AAA_{5}$.
Additionally to \ref {not} we assume that $\rho(X)^G=1$.
By Proposition \ref{Pr1} we may assume also that the singularities of
$X$ are worse than Du Val.

By \cite{DI}, there are three cases: $X_{d+1}\simeq\PP^1$,
$X_{d+1}\simeq\PP^2$ or $X_{d+1}$ is a del Pezzo surface $S_{5}$ of degree $5$.

\begin{lemma}
\label{contractions} Let $V$ be a normal surface and let $C\subset
V$ be a smooth curve such that $(K_V+C)\cdot C<0$. Then $V$ has at
most three singular points on $C$.
\end{lemma}

\begin{proof}
By the adjunction formula \cite{Shokurov-1992-e-ba} we have
\[
(K_V+C)|_C=K_C+\operatorname{Diff}_C,
\]
where $\operatorname{Diff}_C$ is the different, an effective $\QQ$-divisor supported
in singular points of $V$ lying on $C$. Moreover, the coefficients of $\operatorname{Diff}_C$
are $\ge 1/2$. Since, by our conditions $\deg \operatorname{Diff}_C<-\deg K_C\le 2$,
we get that $\operatorname{Diff}_C$ is supported in at most three points.
\end{proof}

\begin{lemma}
\label{A.C1}
For any $i$, the exceptional divisor of
$\psi_{i}$ has at least five connected components.
\end{lemma}

\begin{proof}
Let $E$ be the $\psi_{i}$-exceptional divisor. Since
$\rho(\bar{X}_{i})^G=2$ and $G=\AAA_5$ is a simple group, $E$ is
either connected or the number of connected components of $\psi_{i}$
is $\ge 5$. Assume that $E$ is connected. Since $E$ is a tree of
rational curves, it is irreducible. So, $E\simeq \PP^1$. By Lemma
\ref{fixed-point} the action of $G$ on $E$ is non-trivial. If
$\bar{X}_{i}$ is smooth along $E$, then $E$ is a $(-1)$-curve and
$\psi_{i}(E)$ is a $G$-fixed smooth point. This contradicts Lemma
\ref{fixed-point}. Therefore, $\bar{X}_{i}$ has at least $5$
singular points on $E$. This contradicts Lemma \ref{contractions}.
\end{proof}

\begin{corollary}
\label{A.Co1} If there is a $G$-fixed point on $X_i$ for some $i$,
then there is a fixed point of $G$ on $X_j$ for any $j\le i$.
\end{corollary}

\begin{proof}
Assume that $X_{i}$ has a fixed point of $G$, say $P$. By Lemma
\ref{A.C1} $\psi_{i-1}$ is an isomorphism over $P$. So
$\phi_{i-1}(\psi_{i-1}(P))$ is a fixed point of $G$ on $X_{i-1}$.
\end{proof}

\begin{lemma}
\label{A.Nep1} Suppose that $\rho(X)^G=1$ and $G$ has no fixed
points on $X$. Then $X_{d+1}$ is not a curve.
\end{lemma}

\begin{proof}
Assume that $X_{d+1}\simeq\PP^1$. By Corollary \ref{A.Co1} the group
$G$ has no fixed points on $X_{d}$. Now we choose $\phi_{i}$
according to the construction  \ref{2.11.1}. By Lemma \ref{CB}
the surface $\bar{X}_{d}$ is singular. Since $\AAA_{5}$ is a simple group and
$G$ has no fixed points on $X_{d}$, we see that the exceptional
divisor of $\phi_{d}$ consists of at least five curves $D_{1},\dots,
D_{k}$, where $k\geq5$. Let $f$ be a general fiber of $\psi_{d}$.
Then $$ 0>f\cdot\left(\pi'^*K_{\bar{X}_{d}}+\alpha\sum_{i=1}^k
D_{i}\right)\geq-2+k\alpha,
$$ where $\alpha$ is the codiscrepancy of $D_{i}$. By Lemma
\ref{Discr} and Proposition \ref{Pr1} we see that
$\alpha=\frac{1}{3}$ and $k=5$. By \cite[Theorem 1-5]{Al2} there is
an irreducible non-singular curve $C\in |-3K_{X_{d}}|\neq\emptyset$.
Put $\bar{C}=\phi_{d}^* C-\sum r_{i}D_{i}$, where $r_{i}\geq 0$.
Since $C$ is not a rational curve, we see that $C$ is not a section
of $\psi_{d}$. Then
$$K_{\bar{X}_{d}}+\frac{1}{3}\bar{C}+\sum_{i=1}^k
\left(\alpha+\frac{r_{i}}{3}\right)D_{i}\equiv 0.$$ Hence
$$0=f\cdot\left(K_{\bar{X}_{d}}+\frac{1}{3}\bar{C}+\sum_{i=1}^k
(\alpha+\frac{r_{i}}{3})D_{i}\right)\geq
-2+\frac{2}{3}+\frac{5}{3}=\frac{1}{3},$$ a contradiction.
\end{proof}

\begin{claim}
\label{S5} A smooth del Pezzo surface
of degree $5$ contains exactly five
pencils of conics.
\end{claim}

\begin{proof}
Each pencil of conics $|C|$ has exactly three degenerate members
which are pairs of meeting lines. Since a del Pezzo surface of
degree $5$ contains $15$ such pair of lines, we are done.
\end{proof}

\begin{lemma}
\label{A.Nep2}
Suppose that $\rho(X)^G=1$ and $G$ has no fixed
points on $X$. Then $X\simeq\PP^2$ or $X$ is a del Pezzo
surface $S_{5}$ of degree $5$.
\end{lemma}

\begin{proof}
By Lemma \ref{A.Nep1} and \cite{DI} we have $X_{d+1}\simeq\PP^2$ or
$X_{d+1}\simeq S_{5}$. Now we choose $\phi_{i}$ according to the
construction \ref{2.11.1}. By Proposition \ref{Pr1} we see that
$X_{d}$ has at least one non-Du Val singularity. Hence the
exceptional divisor of $\phi_{d}$ is one orbit $D_{1},\dots, D_{k}$,
where $k\geq 5$. Let $D=\sum_{i=1}^k D_{i}$ and $m:=-D_{i}^2$. Since
the divisor
$$-\phi_{d}^{*}K_{X_{d}}\equiv -K_{\bar{X_{d}}}-\sum_{i=1}^k\alpha
D_{i}$$ is nef and big, so is
\[
-\psi_{d*}\phi_{d}^{*}K_{X_{d}}\equiv
-K_{X_{d+1}}-\sum_{i=1}^k\alpha\psi_{d}(D_{i}).\eqno(*)
\]
By Lemma \ref{Discr}
$\alpha\geq\frac{1}{3}$. Consider two cases:

\begin{case}
$X_{d+1}\simeq\PP^2$. By the above the divisor
$3H-\sum\alpha\psi_{d}(D_{i})$ is nef and big, where $H$ is a line.
Assume that $\deg \psi_{d}(D_{i})\geq 2$. Then $3>2k\alpha$. Since
$k\geq 5$ and $\alpha\geq\frac{1}{3}$, we see that
$2k\alpha\geq\frac{10}{3}$, a contradiction.

Hence, $\psi_{d}(D_{i})$ are lines. Then $k=5$ or $6$. Assume
that $k=5$. Then there is an orbit of five lines on $\PP^2$. Note
that there is an invariant conic $C\subset\PP^2$. The divisor
$\sum\psi_{d}(D_{i})$ meets $C$ in at most $10$ points. Hence there
is an orbit on $C$ consisting of at most $10$ points. However, the
order of any orbit on $C\simeq \PP^1$ is at least $12$, a
contradiction.

Thus $k=6$. Hence the lines $\psi_{d}(D_{i})$ are in general
position, i.e. every line contains five points of intersection.
Therefore, $m\geq 4$ and so $3>k\alpha\geq 3$. Again we have a
contradiction.
\end{case}

\begin{case}
$X_{d+1}$ is  a del Pezzo surface $S_{5}$ of degree $5$.
Assume that $\psi_{d}(D_{i})$ are
$(-1)$-curves. Then $k=10$ and by $(*)$
\[
1=-K_{S_{5}}\cdot
E_{i}>3\alpha-\alpha=2\alpha.
\]
Thus $\alpha<\frac{1}{2}$. It is well
known that on a del Pezzo surface of degree $5$ every $(-1)$-curve
meets three other $(-1)$-curves. Hence, $m\geq 4$, a contradiction.

Assume that $(\psi_{d}(D_{i}))^2\geq 1$. Then
$5=K_{S_{5}}^2>3k\alpha$. Since $k\geq 5$ and
$\alpha\geq\frac{1}{3}$, we see that $3k\alpha\geq 5$, a
contradiction.

Therefore, $(\psi_{d}(D_{i}))^2=0$, i.e. $\psi_{d}(D_{i})$ is a
conic. Then $2k\alpha<5$. Therefore, $m=3$ and $k=5\text{ or }6$. By
Claim \ref{S5} there are only five linear systems of conics. If
$\psi_{d}(D)$ contains two conics of one pencil, then $\psi_{d}(D)$
has at least $10$ components, a contradiction. Therefore,
$\psi_{d}(D)$ consists of five conics contained in different linear
systems. Every component of $\psi_{d}(D)$ meets other four
components. Then $\psi_{d}$ extracts four points on every component
of $\psi_{d}(D)$. Hence, $m=4$, a contradiction.
\end{case}
\end{proof}

\begin{lemma}
\label{A.Fib} Suppose that $\rho(X)^G=1$. Assume that $G$ has a
fixed point $P$ on $X$. Choose $\phi_{0}$ as in the case
\ref{2.11.2}. Assume that $X_{1}=\PP^1$. Then $\bar{X}_{0}$ is
smooth. Moreover $\bar{X}_{0}\simeq\FF_{n}$ and $X\simeq\PP(1,1,n)$.
\end{lemma}

\begin{proof}
In our case $\psi_{0}:\bar{X}_{0}\rightarrow X_{1}=\PP^1$ is a
rational curve fibration. Let $D_{0}$ be a unique exceptional curve
of $\phi_{0}$. Note that $D_{0}$ is contained into the smooth locus
of $\bar{X}_{0}$. Assume that $D_{0}$ is a section. Then there is no
singular fibers. Hence, by Lemma \ref{CB}, we see that
$\bar{X}_{0}\simeq\FF_{n}$. So, we may assume that $D_{0}$ is not a
section of $\psi_{0}$.

Assume that $P$ is not a Du Val singularity. Let $m\geq 2$ be the
degree of the restriction $\psi_{0}|_{D_{0}}:D_{0}\rightarrow\PP^1$.
Then by the Hurwitz formula $-2=-2m+\deg B$, where $B$ is the
ramification divisor. Thus $\deg B=2m-2$. The divisor $B$ is
$G$-invariant. Hence, $\deg B\geq 12$ and $m\geq 7$. Let $f$ be a
general fiber of $\psi_{0}$. Then $$0<-\phi_{0}^{*}K_{X}\cdot
f=-(K_{\bar{X}}+\alpha D_{0})\cdot f=2-m\alpha,$$ where $\alpha$ is
the codiscrepancy of $D_{0}$. Hence, $m<\frac{2}{\alpha}$. By Lemma
\ref{Discr} $\alpha\geq\frac{1}{3}$, a contradiction.

Therefore, $P$ is a Du Val singularity. By Proposition \ref{Pr1} $X$
also has a non-Du Val singular point. Apply construction
\ref{Const}, the case \ref{2.11.1} to $\bar{X}_{0}$ over the base
$\PP^1$: $$ \xymatrix{& \tilde{X}_{0}\ar[dl]_{\xi_0}
\ar[dr]^{\eta_0}&
&& &\tilde{X}_{e}\ar[dl]_{\xi_{e}} \ar[dr]^{\eta_{e}} \\
**[r]\bar{X}_{0} && \bar{X}_{1}&\dots & \bar{X}_{e} &&
\bar{X}_{e+1}.}
$$ Here $\bar{X}_{e+1}$ is smooth and $\bar{X}_{e}$ is singular. We claim that $\bar{X}_{e+1}$ has no
section $M$ with $M^2=-p$, where $p\geq 2$. Indeed, let
$\pi:Y\rightarrow X_{0}$ be the minimal resolution and let $D_{0}$
be a unique exceptional curve over $P$. Note that $\eta_{0}$
contracts curves meeting $D_{0}$. Therefore, $\bar{X}_{e+1}$ has no
invariant $(-p)$-curves $M$ such that $M$ is a section, where $p\geq
2$.

By Lemma \ref{CB} we have $\bar{X}_{e+1}\simeq\PP^1\times\PP^1$. By
Proposition \ref{Pr1} the singularities of $\bar{X}_{e}$ are worse
than Du Val. Let $D_{1},\dots, D_{k}$ be the exceptional curves of
$\xi_{e}$, where $k\geq 5$. Then
$$0<-\left(K_{\tilde{X}_{e}}+\alpha\sum D_{k}\right)\cdot f=2-k\alpha,$$ where
$f$ is a fiber of the projection
$\bar{X}_{e}\simeq\PP^1\times\PP^1\rightarrow\PP^1$ and $\alpha$ is
the codiscrepancy of $D_{i}$. By Lemma \ref{Discr} we have
$\alpha=\frac{1}{3}$, every non-Du Val singularity on $\bar{X}_{e}$
has type $\frac{1}{3}(1,1)$, and $k=5$. According to \cite[Theorem
1-5]{Al2} there is an irreducible non-singular curve $C\in
|-3K_{\bar{X}_{e}}|\neq\emptyset$. Put $\tilde{C}=\xi_{e}^* C-\sum
r_{i}D_{i}$, where $r_{i}\geq 0$. Since $C$ is not rational, $C$ is
not a section. Then
$$K_{\tilde{X}_{e}}+\frac{1}{3}\tilde{C}+\sum
\left(\alpha+\frac{r_{i}}{3}\right)D_{i}\equiv 0.$$ Hence
$$0=f\cdot\left(K_{\tilde{X}_{e}}+\frac{1}{3}\tilde{C}+\sum
\left(\alpha+\frac{r_{i}}{3}\right)D_{i}\right)\geq
-2+\frac{2}{3}+\frac{5}{3}=\frac{1}{3},$$ a contradiction.
\end{proof}

\begin{lemma}
\label{A.Nep3} Suppose that $\rho(X)^G=1$ and $G$ has exactly one
fixed point on $X$. Then either $X\simeq\PP(1,1,2n)$ or
$X\simeq\tilde{\PP^2}_{k,s}$.
\end{lemma}

\begin{proof}
Let $P\in X$ be the fixed point of $G$. Then by Lemma
\ref{fixed-point} $P\in X$ is of type $\frac{1}{r}(1,1)$ for some
$r\geq 2$. Consider two cases:

\begin{case} $r\geq 11$. We choose $\phi_{0}$ as in the
construction  \ref{2.11.2}. Then $X_{1}$ has no fixed points. By
Lemmas \ref{A.Nep2}, \ref{A.Fib} and \ref{CB}, we may assume that
$X_{1}\simeq\PP^2$ or $X_{1}\simeq S_{5}$. Assume that curves of the
exceptional divisor of $\psi_{0}$ contain two singular points or one
non-Du Val singular point. Let $\pi:Y\rightarrow X$ be the minimal
resolution, let $D=\sum D_{i}$ be the exceptional divisor, and let
$D_{0}$ be a unique invariant curve of the exceptional divisor. We
have $D_{0}^2=-r\leq -11$ and there is a morphism
$\pi':Y\rightarrow\bar{X}_{0}$ such that $\pi=\phi_{0}\circ\pi'$.
Hence, by Lemma \ref{Discr2}
$$0<-K_{X}\cdot\pi_{*}E_{i}=-E_{i}\cdot(K_{Y}+D^\sharp)\leq 0,$$
where $E_{i}$ is a $(-1)$-curve contracted by $\psi_{0}\circ\pi'$
and $D^\sharp$ is the codiscrepancy divisor (see \ref{DCoD}), a
contradiction.

Therefore, every contracted curve contains at most one Du Val
singular point of type $A_{p}$. Assume that $X_{1}$ is a del Pezzo surface
of degree $5$. Let
$k:=\rho(\bar{X}_{0}/X_{1})$ and let $l=(\psi_{*}D_{0})^2$. Then
$$
K_{X}^2=10-\rho(Y)+k(p+1)-l-4+\frac{4}{k(p+1)-l}=1-l+\frac{4}{k(p+1)-l}.
$$
So, $l=-1$, $0$, or $1$. On the other hand, $S_{5}$ has no invariant
curves with self-intersection number $-1$, $0$, or $1$, a
contradiction.
Therefore,
$X_{1}\simeq\PP^2$. In this case
$\psi_{*}D_{0}$ is a conic. Hence, $X\simeq\tilde{\PP^2}_{k,s}$.
\end{case}

\begin{case}
$r\leq 10$. Assume that $X$ has a singular fixed point and other
singular points that Du Val singular points. Then we choose
$\phi_{0}$ as in the construction  \ref{2.11.2}.
Then $X_{1}\simeq\PP^2$ or $X_{1}\simeq S_{5}$ and all non-fixed
singular points have type $A_{p}$. Hence,
$$0<K_{X}^2=10-\rho(X_{1})-(p+1)k+m-4+\frac{4}{m},$$ where
$k\in\{12,20,30,60\}$. Therefore, $k=12$, $p=0$, and
$X_{1}\simeq\PP^2$. As above we have $X\simeq\tilde{\PP^2}_{k,0}$.

Now we may assume that $X$ has a fixed singular point and at least
one non-Du Val non-fixed singular point.

Apply construction \ref{Const}. We may construct $\phi_{i}$ as in
the case \ref{2.11.1}. Then $X_{d}$ is a surface with one fixed
point of $G$. Since $r\leq 10$, we see that every $\psi_{i}$ does
not contract curves containing the fixed point of $G$.

Consider the case where $X_{d}\simeq\PP(1,1,2n)$. Assume that there
are non-Du Val singularities on $X_{d-1}$ other than $P$. Let
$D_{1},\dots, D_{k}$ be the exceptional curves of $\phi_{d-1}$. Then
$$
-K_{\PP(1,1,2n)}=(2n+2) l>\alpha\psi_{d-1*}\sum D_{i},
$$
where $l$ is a generator of the Weil divisor class group, and
$\alpha$ is the codiscrepancy of $D_{i}$. Since $\psi_{i}$ does not
contract curves containing the fixed point of $G$, we see that
$\psi_{d}(D_{i})$ also does not contain the fixed point of $G$. Then
$l\cdot\psi_{d-1}(D_{i})\geq 1$. We obtain $n=1$, $m=3$, and
$\alpha=\frac{1}{3}$. Hence, $\psi_{d-1}(D_{i})\geq 2l$, a
contradiction. Therefore, $X_{d-1}$ has exactly one non-Du Val
singularity. Denote it by $P$.

Now we choose $\phi'_{d-1}$ according to the construction
\ref{2.11.2}. By Proposition \ref{Pr1} $X'_{d}$ is a smooth del
Pezzo surface. We obtain $n=1$, $2$, or $3$. Assume that $n=2$ or
$3$. Then
$$0<K_{X_{d-1}}^2\leq 10-\rho(Y)+\frac{8}{3}\leq-15+\frac{8}{3}<0,$$
a contradiction. Hence, $n=1$. Then
$$0<K_{X_{d-1}}^2=
10-\rho(Y)\leq 0,$$ a contradiction.

Consider the case where $X_{d}\neq\PP(1,1,2n)$. Let
$\phi_{d}:\bar{X}_{d}\rightarrow X_{d}$ be the minimal resolution of
$X_{d}$ and let $D$ be the exceptional curve. Let
$\psi_{d}:\bar{X}_{d}\rightarrow X_{d+1}$ be the contraction of
another $G$-equivariant extremal ray. Suppose that $X_{d+1}\simeq
S_{5}$. Then
$$0<K_{X_{d}}^2=10-\rho(\bar{X}_{d})+l-4+\frac{4}{l}.$$ Assume that
$l\geq 4$. Then $(\psi_{d}(D))^2=-1$, $0$, or $1$. On the other hand,
$S_{5}$ has no invariant curves whose self-intersection number
equals to $-1$, $0$, or $1$, a contradiction.

Suppose that $l=3$.
Then
$$
0<K_{X_{d}}^2=10-\rho(\bar{X}_{d})+\frac{1}{3}.
$$
 We see that every
exceptional curve of $\psi_{d}$ meets $D$ in at most two points.
Note that there is some orbit $P_{1},\dots, P_{j}$ consisting of
points of intersections of exceptional curves with $D$. Since the
order of any orbit on $\PP^1$ is at least $12$, we see that the
number of exceptional curves is at least six. Then
$\rho(\bar{X}_{d})\geq 11$, a contradiction.

Assume that $l=2$. Then
$0<K_{X_{d}}^2=10-\rho(\bar{X}_{d})$. We have $\rho(\bar{X}_{d})\geq
5+r$, where $r=\rho(\bar{X}_{d}/X_{d+1})\geq 5$, a contradiction.

Therefore, $X_{d+1}\simeq\PP^2$. Suppose that $l\geq 4$. Then
$\psi_{d}(D)$ is a conic. On the other hand, $\psi_{d}$ is a blowup of at
least twelve points on $\psi_{d}(D)$, a contradiction.

Finally if $\psi_{d}(D)$ is singular, then $\psi_{d}$ is the blowup of
singular points $P_{1},\dots, P_{k}$ of curve $\psi_{d}(D)$. Let
$q:=\mult_{P_{i}}(\psi_{d}(D))\geq 2$ and $\psi_{d}(D)\equiv tH$,
where $H$ is a line on $\PP^2$. Then $D^2=t^2-kq$. By the genus formula,
we have
\[
k\frac{q(q-1)}{2}=\frac{(t-1)(t-2)}{2}.
\]
Hence,
$t^2=kq^2+3t-kq-2$ and so $D^2=kq(q-2)+3t-2>0$, a contradiction.
\end{case}
\end{proof}

\begin{lemma}
\label{A.Nep4a}Assume that $\rho(X)^G=1$ and $X$ has exactly two
singular points. Then $X\simeq F_{2n+k,k-2n,1}$.
\end{lemma}

\begin{proof}
Since $G$ is a simple group, we see that both singular points are
fixed points. By Lemma \ref{fixed-point} there are exactly two fixed
points $P_{1}$ and $P_{2}$ on $X$ and these points have types
$\frac{1}{r_{1}}(1,1)$ and $\frac{1}{r_{2}}(1,1)$ for some $r_{1}$,
$r_{2}$. We may assume that $r_{1}\geq r_{2}$. Apply construction
\ref{Const}. We choose $\phi_{0}$ as in  \ref{2.11.2} with
$P=P_{1}$. By Lemma \ref{A.Nep3} we may assume that
$X_{1}\simeq\PP(1,1,2n)$ or $X_{1}\simeq\tilde{\PP^2}_{k,0}$. Let
$N$ be a unique exceptional curve of $\phi_{0}$. If the exceptional
curves of $\psi_{0}$ contain $P_{2}$, then $r_{2}\geq 5$. Since
$r_{1}\geq r_{2}$, we see that $\psi_{0}(N)$ does not contain the
singular point. Assume that $X_{1}\simeq\tilde{\PP^2}_{k,0}$. Let
$m=-N^2$. Since
$-\phi_{0}^{*}(K_{X}=-K_{\bar{X}_{0}}+\frac{m-2}{m}N$ is nef and
big, we see that $-K_{X_{1}}>\frac{m-2}{m}\psi_{0}(N)$. Since
$\psi_{0}(N)$ does not contain singular point, we see that
$\psi_{0}(N)\cdot E'_{i}\geq 1$, where $E'_{i}$ is the image of the
$(-1)$-curve. Then
$$-K_{X_{1}}\cdot E'_{i}=1-\frac{k-6}{k-4}>\frac{m-2}{m}.$$
On the other hand, $m\geq k-4\geq 8$, a contradiction.

Therefore, $X_{1}\simeq\PP(1,1,2n)$. Since $\psi_{0}$ is the blowup
of one orbit, we see that $N^2=2n-k$, where $k=12$, $20$, $30$, or
$60$. Hence, $X\simeq F_{2n,k-2n,1}$.
\end{proof}

\begin{lemma}
\label{A.Nep4b}Assume that $\rho(X)^G=1$. Then $G$ has at most two
fixed points on $X$.
\end{lemma}

\begin{proof}
Assume that $G$ has three fixed points $P_{1}$, $P_{2}$ and
$P_{3}\in X$. Apply construction \ref{Const} and choose $\phi_{i}$
as in  \ref{2.11.1}. By Lemma \ref{fixed-point} we obtain a surface
$X_{d-1}$ with exactly three singular points. Let $\pi:Y\rightarrow
X_{d-1}$ be the minimal resolution and let $D_{1}$, $D_{2}$, $D_{3}$
be the exceptional curves. Put $n_i:=D_{i}^2$.
\comment{!!!!!!!!!!!!!!!!
We may assume that
$n_{1}\leq n_{2}\leq n_{3}$. Now we choose $\phi_{d-1}$ as in
\ref{2.11.2} with $P=P_{1}$. By Lemma \ref{A.Nep4a}, we see that
$X_{d}\simeq F_{2n,k-2n,1}$. Assume that an exceptional curve
contains a fixed point $P_{j}$. Then each exceptional curve
contains $P_{j}$. Hence, after contraction the order of the local
fundamental group of $P_{j}$ is increased by 12, 20, 30 or 60.
}
By Lemma \ref{fixed-point} all $n_{i}$ are even.
On the other hand, there is a rational
curve fibration $\Phi: Y \rightarrow \PP^1$ such that $D_{1}$,
$D_{2}$, $D_{3}$ are horizontal curves. Assume that $X_{d-1}$ has no
Du Val singularities. Then
$$0<-f\cdot\left(K_{Y}+\frac{n_{1}-2}{n_{1}}D_{1}+
\frac{n_{2}-2}{n_{2}}D_{2}+\frac{n_{3}-2}{n_{3}}D_{3}\right)=-f\cdot\pi^{*}(K_{X_{d-1}}),$$
where $f$ is a generically fiber of $\Phi$. Hence,
$$\frac{n_{1}-2}{n_{1}}+\frac{n_{2}-2}{n_{2}}+\frac{n_{3}-2}{n_{3}}<
2.$$ We obtain $n_{1}=4$ and $n_{2}\leq 6$.

Now apply construction \ref{Const} and we choose $\phi_{d-1}$ as in
\ref{2.11.2} with $P=P_{3}$. We see that the exceptional divisor of
$\psi_{d-1}$ does not contain any singular point. Then $X_{d}$ is a
del Pezzo surface with two fixed points of $G$. Since $n_{1}=4$ and
$n_{2}\leq 6$, we see that $X_{d}$ is not isomorphic to
$F_{2n,k-2n,1}$, a contradiction.

Therefore, $n_{1}=2$. Assume that $n_{2}=2$. Now we choose
$\phi_{d-1}$ as in  \ref{2.11.2} with $P=P_{3}$. We see that the
exceptional divisor of $\psi_{d-1}$ is contained into the smooth
locus. Hence, $X_{d}$ is a del Pezzo surface with two Du Val
singular points, a contradiction with Proposition \ref{Pr1}.
Therefore, $n_{3}\geq n_{2}\geq 4$.

Now we choose $\phi_{d-1}$ as in
\ref{2.11.2} with $P=P_{2}$. We claim that the components of the
exceptional divisor of $\psi_{d-1}$ does not contain singular
points. Indeed, assume that the exceptional divisor of $\psi_{d-1}$
contains $P_{1}$. Then $X_{d}$ is a del Pezzo surface with a fixed
smooth point, a contradiction.

Assume that the exceptional divisor
of $\psi_{d-1}$ contains $P_{3}$. Hence, there is a $(-1)$-curve $E$
on $Y$ meeting both $D_{2}$ and $D_{3}$. Therefore,
\begin{align*}
E \cdot\pi^{*}K_{X_{d-1}}=E\cdot
(K_{Y}+\frac{n_{2}-2}{n_{2}}D_{2}+\frac{n_{3}-2}{n_{3}}D_{3})
\geq\\-1+\frac{n_{2}-2}{n_{2}}+\frac{n_{3}-2}{n_{3}}\geq 0,
\end{align*}
a contradiction.

Thus, if we choose $\phi_{d-1}$ as in \ref{2.11.2} with $P=P_{2}$,
then the exceptional divisor of $\psi_{d-1}$ does not contain the
singular points. The same holds for $P_{3}$. Note that after the
contraction of another $G$-equivariant extremal ray we obtain
$F_{2n,k-2n,1}$. Hence, $n_{2}, n_{3}\in \{ 10, 18, 28, 58\}$. Now
we choose $\phi_{d-1}$ as in \ref{2.11.2} with $P=P_{1}$. Then
$X_{d}\simeq F_{2n,k-2n,1}$. Assume that the exceptional divisor of
$\psi_{d-1}$ contains a singular point. Then $2n=n_{2}-h$ and
$k-2n=n_{3}$ or $2n=n_{2}$ and $k-2n=n_{3}-h$, where $h\in\{ 12, 20,
30, 60 \}$. Since $k\in\{ 12, 20, 30, 60 \}$ and $29\geq n\geq 1$,
this case is impossible. Hence, the exceptional divisor of
$\psi_{d-1}$ does not contain any singular point. We obtain
$F_{2n,k-2n,1}$, where $2n=10$, $18$, $28$, or $58$. Hence,
$n_{2}=n_{3}=10$. Then
$$
K_{X_{d-1}}^2=K_{Y}^2+\frac{(n_{2}-2)^2}{n_{2}}+\frac{(n_{3}-2)^2}{n_{3}}=22-\rho(Y)+\frac{4}{5}.
$$Now we compute $\rho(Y)$. We have $\rho(F_{10,10,1})=20$. Then
$\rho(\bar{X}_{d-1})=20+s$, where
$s=\rho(\bar{X}_{d-1}/F_{10,10,1})\geq 5$ is the number of
components of the exceptional divisor of $\psi_{d-1}$. Then
$\rho(Y)=22+s$. Hence, $K_{X_{d-1}}^2<0$, a contradiction.
\end{proof}

\begin{lemma}
\label{A.Nep4c}Assume that $\rho(X)^G=1$ and $G$ has at least two
fixed points on $X$. Then $X\simeq F_{2n,ak-2n,a}$.
\end{lemma}

\begin{proof}
By Lemma \ref{A.Nep4b} the group $G$ has exactly two fixed points
$P_1$,\, $P_2\in X$ and by Lemma \ref{A.Nep4a} we may assume that $X$ also has a
non-fixed singular point, say $Q$.

First we consider the case where $Q$ is not Du Val. Apply
construction \ref{Const} and choose $\phi_{i}$ as in \ref{2.11.1}.
Then $X_{d}\simeq F_{2n,k-2n,1}$. Let $D_{0,0},\dots, D_{0,l}$ be
the exceptional curves of $\phi_{0}$. Then
$-(K_{\bar{X}_{0}}+\alpha\sum_{j=1}^l D_{j})$ is nef and big, where
$\alpha\geq\frac{1}{3}$ and $l\geq 5$. Let $\tilde{D}_{i,0},\dots,
\tilde{D}_{i,l}$ be the proper transform of $D_{0,0},\dots, D_{0,l}$
on $X_{i}$. Hence, the divisor $-(K_{X_{i}}+\alpha\sum_{j=1}^l\tilde{D}_{i,j})$
is nef and big. Assume that a curve of the exceptional divisor of
$\psi_{i}$ contains a fixed point $P$ and meeting $D_{i-1,j}$. Since
the exceptional curves of $\psi_{i}$ are contained in one orbit, we
see that every curve of the exceptional divisor of $\psi_{i}$
contains a fixed point $P$. Hence, $P$ has type $\frac{1}{m}(1,1)$,
where $m\geq 7$. Let $\pi_{i}: Y_{i}\rightarrow X_{i}$ be the
minimal resolution. Then the divisor
$-(K_{Y}+\frac{m-2}{m}D+\alpha\sum\bar{D}_{i,j})$ is nef and big,
where $D$ is a unique exceptional curve over $P$ and $\bar{D}_{i,j}$
is the proper transform of $\tilde{D}_{i,j}$. On the other hand,
there is $(-1)$-curve $E$ on $Y_{i}$ such that $E$ meets $D$ and
$\bar{D}_{i,j}$. Hence, $0<1-\frac{m-2}{m}-\alpha<0$, a
contradiction. Therefore, $\tilde{D}_{d,0},\dots, \tilde{D}_{d,l}$
do not contain a fixed point. Let $\pi_{d}: Y_{d}\rightarrow
F_{2n,k-2n,1}$ be the minimal resolution. Let $D_{1}$ and $D_{2}$ be
exceptional divisors over fixed points of $G$ and
$\bar{D}_{d,0},\dots, \bar{D}_{d,l}$ be the proper transform of
$\tilde{D}_{d,0},\dots, \tilde{D}_{d,l}$. There is a $G$-equivariant
rational curve fibration $\Phi: Y_{d} \rightarrow \PP^1$ such that
$D_{1}$ and $D_{2}$ are sections. Since curves
$\tilde{D}_{d,0},\dots, \tilde{D}_{d,l}$ do not contain fixed
points, we see that $\bar{D}_{d,0},\dots, \bar{D}_{d,l}$ are
horizontal curves. Hence,
\begin{align*}
0<-f\cdot(K_{Y_{d}}+\frac{2n-2}{2n}D_{1}+\frac{k-2n-2}{k-2n}D_{2}+\alpha\sum_{j=1}^{l}\bar{D}_{d,j})\leq
\\
2-\frac{2n-2}{2n}-\frac{k-2n-2}{k-2n}-l\alpha<0,
\end{align*}
a contradiction.

Therefore, the singularities of $X\setminus \{P_1,\, P_2\}$ are Du
Val. By Lemma \ref{fixed-point} points $P_{1}$ and $P_{2}$ are types
$\frac{1}{m_{1}}(1,1)$ and $\frac{1}{m_{2}}(1,1)$. We may assume
that $m_{1}\geq m_{2}$. If we apply construction \ref{Const} and
choose $\phi_{i}$ as in \ref{2.11.1}. We obtain $X_{d}\simeq
F_{2n,k-2n,1}$. Hence, $m_{1}+m_{2}\geq 12$. Now, we choose
$\phi_{0}$ as in \ref{2.11.2} with $P=P_{1}$. By Lemma \ref{A.Nep3}
we have two possibilities:

\begin{case}
$X_{1}\simeq\tilde{\PP^2}_{k,s}$. Since $m_{1}\geq m_{2}\geq 8$, we
see that the exceptional curves of $\psi_{0}$ do not contain the
fixed point $P_{2}$. Now, we choose $\phi_{1}$ as in \ref{2.11.2}
with $P=P_{2}$. We obtain $X_{2}\simeq\PP^2$. Let $D_{2}$ be a
unique exceptional curve over $P_{2}$. Let $\pi:Y\rightarrow X$ be
the minimal resolution and let $D_{1}$ be the proper transform of a
unique exceptional curve over $P_{1}$. If an exceptional curve of
$\psi_{1}$ meets $D_{1}$, then there is $(-1)$-curves $E$ on $Y$
such that $E$ meets $\bar{D}_{1}$ and $\bar{D}_{2}$. Hence,
$$
0<-E\cdot
\left(K_{Y}+\frac{m_{1}-2}{m_{1}}\bar{D}_{1}+\frac{m_{2}-2}{m_{2}}\bar{D}_{2}\right)
=1-\frac{m_{1}-2}{m_{1}}-\frac{m_{2}-2}{m_{2}}<0,
$$
a contradiction. Therefore, the exceptional curves of $\psi_{1}$ do
not meet $D_{1}$. Then $\psi_{1}(D_{1})$ does not meet
$\psi_{1}(D_{2})$, a contradiction.
\end{case}

\begin{case}
$X_{1}\simeq\PP(1,1,2n)$. As above, we see that the exceptional
curves of $\psi_{0}$ do not contain the fixed point $P_{2}$. Note
that every exceptional curve is a rational curve with one Du Val
singular point of type $A_{s}$. Let
$\pi:\FF_{2n}\rightarrow\PP(1,1,2n)$ be the minimal resolution and
let $D_{2}$ be a unique exceptional curve. Then $D_{2}$ does not
meet $D_{1}$, where $D_{1}$ is the proper transform of a unique
exceptional curve of $\phi_{0}$. Hence, $D_{2}^2=-2n$ and
$D_{1}^2=2n$. Therefore, $X\simeq F_{2n,ak-2n,a}$.
\end{case}
\end{proof}

Now Theorem \ref{1.1} follows from Lemmas \ref{A.Nep2},
\ref{A.Nep3}, and \ref{A.Nep4c}.

\end{document}